\documentstyle[12pt]{article}
\def\ZZ{{{\rm Z}\kern-.52em{\rm Z}}}

\begin{document}

\title{On controllers of prime ideals in group algebras of torsion-free
abelian groups of finite rank}

\author{A.V.TUSHEV}
\date{Department of Mathematics; Dnepropetrovsk National University;
Ul.Naukova 13; Dnepropetrovsk-50; 49050; Ukraine e-mail:
tushev@member.ams.org}
\maketitle
\begin{abstract}
{\it \small
 Let $kA$ be the group algebra of an abelian group $A$ over a field $k$ and let $I$ be an ideal of $kA$.
 We say that a subgroup $B$ of the group $A$ controls the ideal $I$ if $I = (I \cap kB)kA$.
 The intersection $c(I)$ of all subgroups of the group $A$ controlling $I$ is said to be the
  controller of the ideal $I$ . The ideal $I$ is said to be faithful if $I^\dag   = A \cap (1 + I) = 1$.
 In the presented paper we develop some methods which allow us to study controllers of prime faithful
 ideals in group algebras of abelian groups of finite rank. The main idea is that the quotient ring
  $kA/I$ by such an ideal $I$ can be embedded as a domain $k[A]$ in a field $F$ and the group $A$
   becomes a subgroup of the multiplicative group of the field $F$. It allows us to apply for studying
   of $k[A]$ some methods of the theory of fields such as Kummer theory and Dirichlet Unit theorem.
    In its tern properties of $k[A] \cong kA/I$ strongly depend on the properties of the ideal $I$.
    Using these methods, in particular, we obtain an independent proof of a Brookes theorem on controllers
    of prime ideals in the case, where the field $k$ has characteristic
    zero.}

\end{abstract}

      {\bf Introduction}
\vspace{1 pc}\par
      Let $R$ be a ring, let $G$ be a group and let $I$ be a right ideal of the group ring $RG$.
      The ideal $I$ is said to be faithful if $I^\dag   = G \cap (1 + I) = 1$. We say that a subgroup
      $H$ of the group $G$ controls the ideal $I$ if

         $$I = (I \cap kH)kG .~~~~~~~~~~~~~~~~~~~~~~(1)$$

The intersection $c(I)$ of all subgroups of the group $G$
controlling $I$ is said to be the controller of the ideal $I$. It
is easy to note that a subgroup $B \le G$ controls the ideal $I$
if and only if $c(I) \le B$.\par
    Let $H$ be a subgroup of the group $G$ and let $U$ be a right $RH$-module. Since the group ring
    $RG$ can be considered as a left $RH$ -module, we can define the tensor product $U \otimes _{RH} RG$
    (see [5, § 5.1]), which, by [5, § 5.1, proposition 3], is a right $RG$-module named the $RG$-module induced from the $RH$-module $U$. By [2, Chap. III, proposition 5.3], if $M$ is an $RG$
-module and $U \le M$, then

           $$M = U \otimes _{RH} RG~~~~~~~~~~~~~~~~~(2)$$

if and only if

           $$M =  \oplus _{t \in T} Ut,~~~~~~~~~~~~~~~~~~~~~(3)$$

where $T$ - is a right transversal to subgroup $H$ in $G$.\par
      If a $kG$-module  $M$ of some representation  $\varphi $ of a group $G$ over a field $k$ is induced
      from some  $kH$-module $U$, where $H$ is a subgroup of the group $G$, then we say that the representation
      $\varphi $ is induced from a representation of subgroup $H$, where $U$ is the module of the representation
      $\phi $.\par
      Suppose that $M = aRG$ is a cyclic $RG$-module generated by a nonzero element $a \in M$.
      Put $I = Ann_{kG} (a)$ and let $U = akH$, where $H$ is a subgroup of the group $G$. It is not difficult
      to note that, in these denotations, the equation (1) holds if and only if the equation (3) holds. Thus,
      in this case, the equations (1), (2), (3) mean the same.\par
      If $k$ is a subfield of a field $f$ and $A$ is a subset of $f$ then $\left[ {f:k} \right]$ is the dimension
      of the field $f$ over $k$; by $k(A)$ we denote the field generated by $k$ and $A$;  by $k[A]$ we denote the
      domain generated by $k$ and $A$. As usually, $f^* $ denotes the multiplicative group of the field $f$. Let
      $A$ be an abelian group then $t(A)$ denotes the torsion subgroup of  $A$. If the group $A$ is torsion then
      $\pi (A)$ is the set of prime divisors of orders of elements of  $A$.\par
          Let $k$ be a field and let $A$ be an
      abelian group of finite rank. In the presented paper we consider properties of controllers of faithful
      prime ideals in the group algebra $kA$. The main idea of our studying  is that in the case, where $P$ is
      a faithful prime ideal of $kA$, the quotient ring $kA/P$ can be embedded as a domain $k[A]$ in
      a field $F$ and, as the ideal $P$ is faithful, the group $A$ becomes a subgroup of the multiplicative
      group of the field  $F$. It allows us to apply for studying of $k[A]$ some methods of the theory of
      fields such as Kummer theory and Dirichlet Unit theorem. In its tern properties of $k[A] \cong kA/P$
      strongly depend on the properties of the ideal $P$. To prove equation (1) we prove equations (2), (3)
      for a $kA$-module $k[A]$.\par
    In section 1 we consider relations between Kummer theory and equations  (2) and (3) for modules over
    abelian groups. The main result of this section (theorem 1.3) can be considered as a generalization of
    [6, Chap. VIII, theorem 13] to the case of infinite dimensional
    extensions.\par
      In section 2 we study properties of multiplicative groups of certain fields which are generated by a
      finite set and by roots from 1. The most important is that the multiplicative group of such a field
      is a direct product of Chernikov and free abelian groups (see proposition  2.3).\par
      A combination of results of sections 1 and 2 allows us to study controllers of faithful prime ideals
      in group algebras of abelian groups of finite rank. For the first time, such methods were applied in
      [9] (see [9, lemmas 2, 5]) where we proved the equation (3) for modules over abelian groups of finite
      rank. Then the methods were developed in [8, 10, 11, 12]. The main result (theorem 3.1) states that
      in the case where the field $k$ is finitely generated of characteristic zero, the controller of any
      faithful prime ideal $P$ of the group algebra $kA$ of an abelian group $A$ of finite rank with
      Chernikov torsion subgroup is finitely generated. In the case, where the group $A$ is minimax, such
      a result follows from a Segal theorem [8, theorem 1.1]. An abelian group is said to be minimax
      if it has a finite series each of whose factor is either cyclic or   quasi-cyclic.\par
     In section 4 we obtain an independent proof of a Brookes theorem [1, theorem A] in the case of the
     field of characteristic zero (theorem 4.4). As it became known recently, the original proof of
     [1, theorem A] is incorrect. Moreover, there is an example which shows that in fact the theorem
     does not hold in the case of the field of positive characteristic. So, a new independent proof
     of the Brookes theorem in the case of the field of characteristic zero is quite topical.
\vspace{1 pc}

     {\bf 1. Kummer theory and induced modules.}
\vspace{1 pc}\par
     Let $k$ be a subfield of a field $f$ and let $G$ be a subgroup of the multiplicative group $f^*$ of the
     field $f$. Then the field $k(G)$ may be considered as a $kG$-module and the field $k$ can be considered
     as a $k(G \cap k^*)$-module. Therefore, we can define the tensor product $k \otimes _{k(G \cap k^*)} kG$
     and the equation $k(G) = k \otimes _{k(G \cap k^*)} kG$ means that $k(G) =  \oplus _{t \in T} kt$, where
     $T$ is a transversal to the subgroup $G \cap k^*$ in the group $G$. If
     $\left| {G/G \cap k^*} \right| = m < \infty $ the equation $k(G) = k \otimes _{k(G \cap k^*)} kG$ holds
     if and only if $\left[ {k(G):k} \right] = m$. The relations between $\left| {G/G \cap k^*} \right|$ and
     $\left[ {k(G):k} \right]$ were considered in Kummer theory (see [6, Chap. VIII, theorem 13]). So, we
     can hope that there are relations  between induced modules over abelian groups and Kummer theory. These
     relations are studied in this section.
\vspace{1 pc}\par
     {\bf Lemma 1.1.} {\it Let $k$ be a subfield of a field $f$ and let $G$ be a subgroup of $f^*$ such that $k^* \le G$ ,
     $k$ contains a primitive root from 1 of degree 4 if the quotient group $G/k^*$ has an element of order 4 ,
     the quotient group $G/k^*$ is torsion, for any $p \in (\pi (t(G)) \cap \pi (G/k^*))$ the field  $k$ contains
     a primitive root from 1 of degree $p$   and $char k \notin \pi (G/k^*)$. Let $g \in G\backslash k^*$ and let
     $\bar g$ be the image of $g$ in the quotient group $G/k^*$. Suppose that $\left| {\bar g} \right| = t$ then
     $\left[ {k(g):k} \right] = t$.}\par
      {\bf Proof.} Let $g^t  = a \in k^* $. Suppose that $a = 1$ then $g \in t(G)$ and hence all prime divisors of
      $t$ belong to the set $\pi (t(G)) \cap \pi (G/k^*)$. But then for any prime divisor $p$ of $t$ the field
      $k$ contains a primitive root from 1 of degree $p$, and it easily implies that $a \ne
      1$.\par
      Suppose that for some prime divisor $p$ of $t$ there is an element $b \in k$, such that $b^p  = a$. Then
      $(g^m b^{ - 1} )^p  = 1$, where $m = t/p$. Evidently, $p \in (\pi (t(G)) \cap \pi (G/k^*))$ and hence $k$
      contains a primitive root from 1 of degree $p$. Therefore, $g^m b^{ - 1}  \in k^* $ and, as $b \in k^* $,
      we have $g^m  \in k^* $. But it is impossible because   $\left| {\bar g} \right| =
      t$.\par
      Suppose that 4 divides  $t$ and $a \in  - 4k^4 $ then there is an element $b \in k$ such that
      $a =  - 4b^4 $      and the quotient group $G/k^*$ has an element of order 4. As the quotient group $G/k^*$
      has an element of
      order 4, the field $k$ contains  a primitive root from 1 of degree 4. Let $m = t/4$ and $h = g^m $
      then $h^4  =  - 4b^4 $ and hence $h^2  = ( \pm 2i)b^2 $, where $i$ is a primitive root from 1 of degree 4.
      Since the field $k$ contains  $i$, we have $h^2  = g^{2m}  \in k^* $. But it is impossible because
      $\left| {\bar g} \right| = t$.\par
      Thus, $a \notin k^p $ for any prime divisor $p$ of $t$ and $a \notin  - 4k^k $ if  $4$ divides $t$. Then it
      follows from [6, Chap. VIII, theorem 16] that the element $g$ is a root of an irreducible polynomial
      $X^t  - a$ over the field $k$ and, by [6, Chap. VII, proposition 3], $[d:k] = t$.\vspace{1 pc}\par

     {\bf Lemma 1.2.} {\it Let $k$ be a subfield of a field $f$ and let $G$ be a subgroup of $f^*$ such that $k^* \le G$,
      $k$ contains a primitive root from 1 of degree 4 if the quotient group $G/k^*$ has an element of order 4,
      the quotient group $G/k^*$ is torsion,  for any $p \in (\pi (t(G)) \cap \pi (G/k^*))$ the field  $k$
      contains a primitive root from 1 of degree $p$   and $char k \notin \pi (G/k^*)$. Suppose that
      $\left| {G/k^*} \right| = n$, where $n = p^m $ and $p$ is a prime number. Then
      $\left[ {k(G):k} \right] = p^m  = n$ .}\par
      {\bf Proof.} Suppose that the field $k$ contains a primitive root $\xi $ from 1 of degree $p$. The proof is
       by induction on $m$.\par
        Let $g$ be an element of $G$ such that $\left| {\bar g} \right| = p$, where $\bar g$ is the image of
        $g$ in the quotient group $G/k^*$. Let $d = k(g)$ then, by lemma 1.1, $[d:k] = p$. Suppose that there
        is an element $h \in G \cap d^*$ such that  $\left| {\bar h} \right| = p^2 $, where $\bar h$ is the
        image of $h$ in the quotient group $G/k^*$.\par
        By lemma 1.1, $[k(h):k] = p^2 $ and, as $k(h) \le k(g)$ and $[k(g):k] = p$, we obtain a contradiction
        with [6, Chap. VII, proposition 2]. Thus, we can assume that the quotient group $(d^*  \cap G)/k^* $
        is elementary abelian. Put $d^*  \cap G = L$, then $k^*  \le L$ and , evidently,
        $d = k(L)$. By [6, Chap. VIII, theorem 13], $p = [d:k] = |L^p /(k^*)^p |$ and, as the field $k$
        contains a primitive root from $1$ of degree $p$, it is not difficult to show that $|L/k^*| = p$. Then
        $G \cap d^* = k^* < g > $ and hence  $\left| {G/d^*} \right| = p^{m - 1} $. Therefore, by the induction
        hypothesis, $\left[ {d(G):d} \right] = p^{m - 1} $. Evidently, $d(G) = k(G)$ and, as
        $\left[ {d:k} \right] = p$, it follows from [6, Chap. VII, Proposition 2] that
        $\left[ {k(G):k} \right] = p^m  = n$.\par
      Suppose that the field $k$ does not contain $\xi $. Put $d = k(\xi )$ then $[d:k] = r < p$
      because the cyclotomic polynomial for $\xi $ has degree $p-1$.\par
      Now, we show that $G \cap d^* = k^*$. Suppose that there is an element $g \in (G \cap d^*)\backslash k^*$.
      Evidently, there is no harm in assuming that $\left| {\bar g} \right| = p$, where $\bar g$ is the
      image of $g$ in the quotient group $(G \cap d^*)/k^*$. Then, by lemma 1.1, $[k(g):k] = p$ but, as
      $k(g) \le d$, we obtain a contradiction with $[d:k] = r < p$. Thus, $G \cap d^* = k^*$ and hence
      $Gd^*/d^* \cong G/(G \cap d^*) = G/k^*$. Therefore, $|G/(d^*\cap G)| = n = p^m $. Then, as it was
      shown above, $[d(G):d] = p^m $ and hence, by [6, Chap. VII, Proposition 2], $[d(G):k] = p^m r$,
      where $[d:k] = r < p$. Since $d(G) = k(G)(\xi )$, we have $[d(G):k(G)] \le r$. Then it follows
      from [6, Chap. VII, proposition 2] and the equation $[d(G):k] = p^m r$ that $[k(G):k] \ge n$.
      But the relation $[k(G):k] > n$ is impossible, because $[k(G):k] \le \left| {G/k^*} \right| = n$.
      Thus, $[k(G):k] = n$.\vspace{1 pc}\par

      {\bf Theorem 1.3.}{\it Let $k$ be a subfield of a field $f$ and let $G$ be a subgroup of $f^*$ such that
      $k^* \le G$, $k$ contains a primitive root from 1 of degree 4 if the quotient group $G/k^*$ has
      an element of order 4, the quotient group $G/k^*$ is torsion, for any $p \in (\pi (t(G)) \cap \pi (G/k^*))$
      the field  $k$ contains a primitive root from 1 of degree $p$   and $char k \notin \pi (G/k^*)$.  Then
      $k(G) = k \otimes _{k(k^*)} kG$ =$ \oplus _{t \in T} kt$, where $T$ is a transversal to $k^*$ In
      $G$.}\par
      {\bf Proof.} Since the quotient group $G/k^*$ is locally finite, it is sufficient to show that for
      any subgroup
      $H \le G$ such that $k^* \le H$ and $\left| {H/k^*} \right| < \infty $ the equation
      $k(H) = k \otimes _{k(k^*)} kH$ holds. Thus, there is no harm in assuming that
      $\left| {G/k^*} \right| = n < \infty $ and $char k$ does not divide $n$.
      Since the equation $k(G) = k \otimes _{k(k^*)} kG$ means that $k(G) =  \oplus _{t \in T} kt$,
      where $T$ is a transversal to $k^*$ in the group $G$, it is sufficient to show that $[k(G):k] = n$.
      We use the induction on $n$.\par
      Let $p$ be the smallest prime divisor of $n$ and let $N/k^*$ be the
      Sylow $p$-subgroup of the quotient group $G/k^*$. Suppose that $|N/k^*| = t$. Let $d = k(N)$,
      then by lemma 1.2, $[d:k] = t$.\par
      Now, we show that $d^* \cap G = N$. Suppose, that there is an element $g \in (d^* \cap G)\backslash
      N$. Evidently, there is no harm in assuming that $\left| {\bar g} \right| = p'$, where $\bar g$ is
      the image of the element $g$ in the quotient group $G/(k^* \cap G)$ and $p'$ is a prime number
      such that $p < p'$. Then, by lemma 1.1, $[k(g):k] = p'$ and, as $k(g) \le d$, by
      [6, Chap. VII, proposition 2], $p'$ divides $t$ but it is impossible because $p \ne p'$. Thus,
      $d^* \cap G = N$ and hence  $Gd^*/d^* \cong G/(d^* \cap G) = G/N$. So, we can conclude that
      $\left| {Gd^*/d^*} \right| = \left| {G/N} \right| = n/t$. Then by the induction hypothesis,
      $\left[ {d(G):d} \right] = n/t$ and, as $[d:k] = t$, by [6, Chap. VII, proposition 2],
      we have $\left[ {d(G):k} \right] = n$. Finally, it is easy to see that $d(G) = k(G)$
      and the assertion follows.\vspace{1 pc}\par

      {\bf Corollary 1.4.}{\it Let $k$ be a subfield of a field $f$ and let $G$ be a subgroup of $f^*$ such that $k$
      contains a primitive root from 1 of degree 4 if the quotient group $Gk^*/k^*$ has an element of order 4,
      the quotient group $Gk^*/k^*$ is torsion, for any $p \in (\pi (t(Gk^*)) \cap \pi (Gk^*/k^*))$ the field  $k$
      contains a primitive root from 1 of degree $p$   and $char k \notin \pi (G/k^*)$.  Then
      $k(G) = k \otimes _{k(k^* \cap G)} kG= \oplus _{t \in T} kt$, where $T$ is a transversal to
      $k^* \cap G$ in $G$.}\par
      {\bf Proof.}  To prove the corollary, it is sufficient to apply the above theorem to the field $k$ and the
      group $Gk^*$.\vspace{1 pc}\par

     {\bf 2. On multiplicative groups of certain fields.}
\vspace{1 pc}\par

     We will say that a field $k$ is regular if it is countable and the multiplicative group of the field
     $k$ is a direct product of a torsion group and a free abelian group.
\vspace{1 pc}\par
    { \bf Lemma 2.1.}{\it Let $f = k(S)$ be a transcendent extension of a field $k$, where $S$ is
    a finite set of elements of the field $f$. If  the field $k$ is regular then so is the field
    $f$.}\par
     {\bf Proof. }Since the set $S$ is finite, so is the transcendent degree of the field $f$ over the
     subfield $k$. Let $z \notin k$ and let $z,z_1 ,...z_n $ be a maximal system of algebraically
     independent over $k$ elements of the field $f$ . Put $f_1  = k(z,z_1 ,....,z_n )$ then it is
     not difficult to show that $f_1^*  = k^*  \times N$, where $N$ is a countable free abelian
     group and hence, as the field $k$ is regular, so is the field $f_1 $. As $f \ge f_1  \ge k$
     and $f = k(S)$, we can conclude that $f = f_1 (S)$. Then, as the field $f$ is an algebraic
     extension of the field $f_1 $, we see that $\left[ {f:f_1 } \right]$ is finite and hence we
     can put $\left[ {f:f_1 } \right] = m < \infty $.\par
     Let $\varphi _z :f^*  \to f_1^* $, be a homomorphism which maps each element of the multiplicative
     group $f^ *  $ to its regular norm over the field $f_1 $. Then, evidently,
     $\varphi _z (z) = z^m $.\par
     Suppose that $z \in t(f^ *  /Ker\varphi _z )$ then there is a positive integer $r$ such that
     $z^r  \in Ker\varphi _z $ and hence $1 = \varphi _z (z^r ) = (\varphi _z (z))^r  = z^{mr} $.
     Therefore, $z \in t(f^ *  )$. On the other hand, as the field $f$ is a transcendent extension
     of the field $k$, the quotient group $f^* /k^* $ is torsion-free and , as $z \notin k$, we see
     that $z \notin t(f^* )$. Thus, a contradiction is obtained and hence $z \notin t(f^ *  /Ker\varphi _z )$.
     Since $f^ *  /Ker\varphi _z  \cong \varphi _z (f^ *  ) \le f_1^ *  $ and the field $f_1 $ is regular,
     it is not difficult to show that the quotient group $f^ *  /f_x $ is free abelian , where $f_x $ is a
     subgroup of  $f^ *  $ and $f_x /Ker\varphi _z  = t(f^ *  /Ker\varphi _z
     )$.\par
     Thus, for any element $z \notin f\backslash k$ there exists a subgroup $f_z $ of the group $f^* $
     such that $z \notin f_z $ and the quotient group $f^* /f_z $ is a free abelian. Put
     $T = \bigcap\limits_{z \in f\backslash k} {f_z }$, evidently $T \le k^* $ and the quotient group
     $f^* /T$ is a residually free abelian and, as the group $f^* $ is countable, it implies that the
     quotient group $f^* /T$ is free abelian. Since $T \le k^* $, we see that $T$ is a direct product
     of a torsion and a free abelian groups and hence the field $f$ is regular.\vspace{1 pc}\par

     Let $k$ be a field, the quotient group $\bar k = k^*/t(k^*)$ is said to be the reduced
     multiplicative group of the field $k$.\vspace{1 pc}\par

     {\bf Lemma 2.2.}{\it Let $f$ be a finite extension of a field $k$ such that $\left[ {f:k} \right] = n$.
     Suppose that the field $f$ is regular. Then  $\bar f = F \times (\bar f \cap \bar k^{1/n} )$, where $F$
     is a countable free abelian group.}\par
     {\bf Proof.}  Let $\varphi :f^ *   \to k^ *$ be a homomorphism given by $\varphi :x \mapsto \left| x \right|$,
     where $\left| x \right|$ is the regular norm of an element $x \in f$ over $k$. It is easy to note that
     $\varphi $ induces a homomorphism $\bar \varphi :\bar f \to \bar k$. It follows from the definition of
     the homomorphism $\bar \varphi $ that $\bar \varphi (x) = x^n $ for any $x \in \bar k$. Since the field $f$
     is regular, we see that the groups $\bar k$ and $Ker\bar \varphi $ are free abelian. Put
     $A = Ker\bar \varphi $, since $\bar \varphi (x) = x^n $ for any $\bar x \in k$, we can conclude that
     $1 = A \cap \bar k$ and hence $1 = A \cap (\bar f \cap \bar k^{1/n})$.\par
     Let $x \in \bar f$ and $\bar \varphi (x) = y \in \bar k$ then
     $\bar \varphi (x^n y^{ - 1} ) = \bar \varphi (x^n )\bar \varphi (y^{ - 1} ) =
      \bar \varphi (x)^n (y^{ - 1} )^n  = y^n (y^{ - 1} )^n  = 1$.  Thus, for any element $x \in \bar f$
      there is an element $y \in \bar k$ such that $x^n y^{ - 1}  \in A$. Therefore,
      $(\bar f)^n  \le \bar k \times A$ and hence the quotient group $R = \bar f/(\bar f \cap \bar k^{1/n})$
      is torsion free. Since $1 = A \cap (\bar f \cap \bar k^{1/n} )$ and the subgroup $A$ is free abelian,
      we can conclude that the quotient group $R$ has a free abelian subgroup
      $T = (A(\bar f \cap \bar k^{1/n} ))/(\bar f \cap \bar k^{1/n} )$ such that $R^n  \le T$.
      Then, as the group $R$ is torsion free, it easily implies that the group $R$ is free abelian.
      Thus, the quotient group $R = \bar f/(\bar f \cap \bar k^{1/n} )$ is free abelian and hence
      there is a free abelian subgroup $F \le \bar f$ such that $\bar f = F \times (\bar f \cap \bar k^{1/n} )$.
\vspace{1 pc}\par

     {\bf Proposition 2.3.}{\it Let $f$ be a field of characteristic zero generated by a finite set $S$
     and by all roots from 1 of degree $p^n $, where $n$ are positive integers, $p \in \pi $ and $\pi $
     is a finite set of prime numbers. Then  the set $\pi (t(f^*))$ is finite  and  the field $f$ is
     regular, $f^* = T \times F$, where   $F$is a free abelian  group and $T$:\\
     (i)  is a finite group if  the field $f$ is finitely generated (that is if $\pi = \emptyset$);\\
     (ii) is a locally cyclic Chernikov group.}\par
     {\bf Proof.}  Let $k$ be a subfield of $f$ generated by all roots from 1 of degree $p^n $, where $n$
     are positive integer and $p \in \pi $, then $f = k(S)$. If $\pi  = \emptyset $ then $k$ is the
     field of rational numbers. Let $k_1 $ be the maximal algebraic extension of $k$ in $f$. Since
     $\left| S \right| < \infty $, we see that $\left[ {k_1 :k} \right] < \infty $ and $f$ is a
     finitely generated transcendent extension of $k_1 $. As $f$ is a finitely generated transcendent
     extension of $k_1 $, we see that $t(f^*) = t(k_1 ^*)$ and it follows from lemma 2.1 that it is
     sufficient to show that the field $k_1 $ is regular. So, we can assume that the elements of
     $S$ are algebraic over $k$ and hence $\left[ {f:k} \right] = r < \infty $.\par
     (i) In this case, $f$ is a finitely generated algebraic field of characteristic zero and hence
     $f$ is a finite extension of the rational number field. Then, be the Dirichlet Unit theorem
     (see [4, Chap. I, theorem 13.12]), the group $U$ of  unites of the ring $R$ of algebraic
     integers of the field $f$ is finitely generated. There exists a homomorphism $\varphi $
     of the group $f^*$ into the free abelian group $I(R)$ of  fractional ideals of  $R$ which
     maps each element $a \in f^*$ into the fractional ideal generated by the element $a$
     (see [4, Chap. I, section 4]) and $Ker\varphi  = U$. Then, by the theorem on homomorphism,
     the quotient group $f^*/U$ is free abelian and, as the group $U$ is finitely generated, it
     easily implies that  $f^* = T \times F$, where   $F$ is a free abelian  group and $T$   is
     a finite group.\par
     (ii) At first, we show that the set $\pi (t(f^*))$ is finite and hence the group $t(f^*)$ is
     Chernikov. Suppose that the field $f$ contains a primitive root $\xi $ from 1 of degree
     $q \notin \pi $ . It is well known that the cyclotomic polynomial for $\xi $ is irreducible
     over $k$ and hence $\left[ {k(\xi ):k} \right] = q - 1$. Then $q - 1 \le r$ and evidently the
     set of all such $q$ is finite and, as the set $\pi $ is finite, so is $\pi (t(f^*))$. Therefore,
     $T = t(f^*)$ is a Chernikov locally cyclic group.\par
     Now, we show that the field $f$ is regular. Let $D$ be the set of all roots from 1 of degree
     $p^2 $, where $p \in \pi $, then the set $D$ is finite. Let $h$ be a subfield of $f$ generated
     by $(S \cup D)$. Then, by (i),  $h^*$ is a direct product of a finite and a free abelian groups.
     Evidently, the field $f$ has an infinite series $\left\{ {h_i| i\in I } \right\}$ of subfields
     such that $h = h_1 $, $h_i  \le h_{i + 1} $, $ \cup _{i \in I } h_i  = f$ and
     $h_{i + 1}  = h_i (\zeta _i )$, where $\zeta _i $ is a root from 1 of degree $p^n $
     for some $p \in \pi $ and positive integer $n$. We also can assume that
     $\zeta _i^p  \in h_i $ then it follows from lemma 1.1 that $\left[ {h_{i + 1} :h_i } \right] = p$ if
     $h_i  \ne h_{i + 1} $. Evidently, $\bar h_i  \le \bar h_{i + 1} $, $ \cup _{i \in I} \bar h_i  = \bar f$
     and it follows from (i) that each group $\bar h_i $ is  free abelian. Then it is sufficient to show
     that the quotient group $\bar h_{i + 1} /\bar h_i $ is free abelian for each
     $i$.\par
     By  lemma 2.2, $\bar h_{i + 1}  = F_i  \times (\bar h_{i + 1}  \cap \bar h_i ^{1/p} )$, where $F_i $
     is a countable free abelian group and it is sufficient to show that
     $\bar h_i  = \bar h_{i + 1}  \cap \bar h_i ^{1/p} $. Suppose that
     $\bar h_i  \ne \bar h_{i + 1}  \cap \bar h_i ^{1/p} $ then there is an element
     $a \in \bar h_{i + 1} \backslash \bar h_i $ such that  $a^p  \in \bar h_i $
     and hence there is an element $b \in h_{i + 1} \backslash h_i $ such that
     $b^p  \in h_i $ and $b \notin \left\langle {\zeta _i } \right\rangle $. It
     easily implies that the group $h_{i + 1} ^*$ has a subgroup $G$ such that
     $h_i ^* \le G$ and $G/h_i ^*$ is an elementary abelian $p$-group of order $p^2 $.
     Therefore, as $h$ contains a primitive root from 1 of degree $p^2 $,  it follows
     from lemma 1.2 that $\left[ {h_i (G):h_i } \right] = p^2 $ but it is impossible because
     $h_i ^*(G) \le h_{i + 1} ^*$ and $\left[ {h_{i + 1} :h_i } \right] =
     p$.\par
    Thus, the field $f$ is regular and hence $f^* = T \times F$, where $T$ is a locally cyclic
    Chernikov group and $F$ is a free abelian group.\vspace{1 pc}\par

{\bf 3. Controllers in prime ideals of abelian groups of finite}
rank. \vspace{1 pc}\par

      Let $A$ be an abelian group and let $B$ be a subgroup of $A$. The set $is_A (B)$ of elements
      $a \in A$ such that $a^n  \in B$, fore some positive integer $n$, is a subgroup of $A$ which
      is said to be the isolator of the subgroup $B$ in the group $A$. The subgroup $B$ is said to
      be dense if  $is_A (B) = A$ and the subgroup $B$ is said to be isolated if $is_A (B) = B$.
\vspace{1 pc}\par
      {\bf Theorem 3.1.}{\it Let $k$ be a finitely generated field of characteristic zero, let $G$ be an
      abelian group of finite rank such that the torsion subgroup $t(G)$ is Chernikov and let
      $P$ be a prime faithful ideal of $kG$. Then the controller of $P$ is finitely
      generated.}\par
      {\bf Proof.}  Let $M$ be the field of fraction the domain  $kG/P$, then $M = k(G)$ and $G$ is a
      subgroup of the multiplicative group of $M$. Let $h$ be the algebraic closure of the field $M$.
      Let $d = k(H,i)$, where $H$ is a finitely generated dense subgroup of $G$ which contains elements
      of order $p$ for each $p \in \tau  = \pi (t(G))$, and $i$ is a primitive root of degree 4 from 1.
      By proposition 2.3(i), $d^* = T \times L$, where $T$ is a finite group and $L$ is a free abelian group.
      It implies, that $G \cap d^*$ is a finitely generated dense subgroup of $G$ and changing $H$ by
      $G \cap d^*$ we can assume that $H = G \cap d^*$.\par
      Let $D = Gd^*$ then $D/d^* = Gd^*/d^* = G/(G \cap d^*) = G/H$ and hence the quotient group $D/d^*$ is
      torsion, besides $d$ contains a primitive root of degree 4 from 1. So, it would be possible to
      apply theorem 1.3 but there may be a situation that not for any $p \in (\pi (t(D)) \cap \pi (D/d^*))$
      the field  $d$ contains a primitive root from 1 of degree $p$. However, as  $D = Gd^*$, where
      $\pi (t(G)) \subseteq \pi (t(d^*))$ and $d^*$ is an almost free abelian group, the above situation
      may be possible only for $p \in \pi  = \pi ((is_{\bar d} \bar H)/\bar H)$, where $\bar d = d^*/t(d^*)$
      and $\bar H = Ht(d^*)/t(d^*)$. Since the group $\bar d$ is free abelian, it is easy to note that the
      set $\pi $ is finite. Let $p$ be the biggest prime from $\pi $ and let $\omega $ be the set of all
      primes $q \le p$.  Let $\xi $ be a primitive root from 1 of degree  $p$, let $d_1  = d(\xi )$
      and let $H_1  = G \cap d_1 ^*$. By lemma 2.2, $\bar d_1  = R \times (\bar d_1  \cap \bar d^{1/n} )$,
      where $R$ is a free abelian group and $n = \left[ {d_1 :d} \right] < p$, and hence
      $\pi ((is_{\bar d_1 } \bar H_1 )/\bar H_1 ) \subseteq \varpi \backslash \{ p\} $, where
      $\bar d_1  = d_1 ^*/t(d_1 ^*)$ and $\bar H_1  = (G \cap d_1 ^*)t(d_1 ^*)/t(d_1 ^*)$. Thus,
      after several steps, adding primitive roots from 1, we obtain a field $f$ such that the
      quotient group $F/f^*$ is torsion, where $F = Gf^*$, $i \in f$ and for any
      $p \in (\pi (t(F)) \cap \pi (F/f^*))$ the field  $f$ contains a primitive root from 1 of degree $p$.
      Therefore, by theorem 1.3, $f(F) = \oplus _{t \in T} ft$, where $T$ is a transversal to
      $f^*$ in $F$. Since $F = Gf^*$ we see that $f(G) = \oplus _{t \in T} ft$, and $T$ can be
      chosen as a transversal to $L = f^* \cap G$ in $G$. Therefore, $k(G) = \oplus _{t \in T} k(L)t$
      and hence $k[G] = \oplus _{t \in T} k[L]t$, where $T$ is a transversal to $L$ in $G$.
      Evidently, it implies that $P = (P \cap kL)kG$ and, as the field $f$ is finitely
      generated, it follows from proposition 1.3(i) that so is the subgroup $L$.
\vspace{1 pc}\par

      {\bf Corollary 3.2.} {\it Let $k$ be a finitely generated field of characteristic zero and
      let $G$ be an
      abelian group of finite rank such that the torsion subgroup $t(G)$ is Chernikov. Then any
      faithful irreducible representation of  the group $G$ over the field $k$ is induced from
      some finitely generated subgroup of $G$.}\par
      {\bf Proof.}  Let $M$ be a module of a faithful irreducible representation of  the group $G$ over
      the field $k$. Then $M$ is a simple $kG$-module and hence $M \cong kG/P$, where
      $P = Ann_{kG} (a)$ for some nonzero element $a \in M$ is a faithful maximal ideal of
      $kG$. By the above theorem, the ideal $P$ is controlled by a finitely generated subgroup
      $H \le G$. It means that $M =  \oplus _{t \in T} Ut$, where $U = kH/Ann_{kH} (a)$ and
      $T$ is a transversal to $H$ in $G$,  and hence $M = U \otimes _{kH} kG$ .
\vspace{1 pc}\par
     {\bf Theorem 3.3.}{\it Let $k$ be a finitely generated field of characteristic zero, let $A$ be a
     torsion-free abelian minimax group and let $P$ be a faithful prime ideal of  the group
     algebra $kA$. Let $R = kA/P$ and let $h$ be the field of fractions of $R$ then
     $A \le h^ *  $. Let $C$ be a finitely generated dense subgroup of $A$ and let
     $\tau $ be the set of  all roots from 1 of degree $p^n $ for all $p \in \pi (A/C)$
     and all positive integer $n$ in addition with a primitive root of degree 4 from 1 .
     Let $f$ be a field obtained by addition of all roots from the set $\tau $ to the field
     $h$ and let $s$ be a subfield of $f$ generated by all roots from $\tau $   and by the
     subgroup $C$ and let $B = A \cap s^ *  $. Let $D$ be a dense subgroup of $C$,  let $d$
     be  a subfield of $f$   generated by $D$ and by $t(s^*)$ and let $H = A \cap d^ *  $. Then
     :\\
     (i) $B$ and $H$ are finitely generated dense subgroups of $A$, besides $B \ge
     H$;\\
     (ii) $B$ and $H$ control the ideal $P$.}\par
     {\bf Proof.} (i) By proposition 2.3(ii), the field $s$ is regular and hence so is its subfield $d$.
     It easily implies that $B$ and $H$ are finitely generated dense subgroups of $A$ and the
     relation $B \ge H$ is evident.\par
     (ii) Put $F = As^*$. Since $B$ is a dense subgroup of $A$, we see that the quotient group
     $F/s^*$ is torsion. By the definition of $\tau $, $s$ contains a primitive root from 1 of
     degree 4 and a primitive root of degree $p$ for any $p \in \pi (F/s^*)$. So, we can conclude
     that $F$ and $s$ meet all conditions of theorem 1.3 and hence $s(A) = s \otimes _{s(s^*)} sA$.
     Therefore, $k[A] = k[B] \otimes _{kB} kA$ and hence $P = (P \cap
     kB)kA$.\par
      Put $L = Ad^*$. Since $H$ is a dense subgroup of $A$, we see that $L/d^*$ is a torsion group.
      As $d \le s$ and $t(s^*) \le d$, we can conclude that $t(s^*) = t(d^*)$ and hence $d$ contains
      a primitive root from 1 of degree 4.\par
      Now, we show that for any $p \in (\pi (t(L)) \cap \pi (L/d^*))$
      the field  $d$ contains a primitive root from 1 of degree $p$. Evidently,
      $L/d^* = Ad^*/d^* \cong A/(A \cap d^*) = A/H$ and the group $A$ has a series of subgroups
      $A \ge B \ge H$. Since $\pi (A/B) \subseteq \tau $ and $\tau  \subseteq \pi (t(s^*)) = \pi (t(d^*))$,
      it is sufficient to consider the case where $p \in (\pi (t(Bd^*)) \cap \pi (B/H))$. But $Bd^* \le s^*$
      and, as $t(s^*) = t(d^*)$, we can conclude that $p \in \pi (t(d^*))$. Thus, $L$ and $d$ meet all
      conditions of theorem 1.3 and the above arguments show that the subgroup $H$ controls the ideal $P$.
\vspace{1 pc}\par
     {\bf 4. On controllers of standardized prime faithful ideals in group algebras of abelian
     groups of finite rank.}
\vspace{1 pc}\par
     Let $k$ be a field, let $A$ be a torsion-free abelian group of finite rank acted by a group $\Gamma $
     and let $I$ be an ideal of a group algebra $kA$. The set $\Delta _\Gamma  (A)$ of elements of $A$
     which have finite orbits under action of the group $\Gamma $ is a subgroup of the group $A$. The
     subgroup $N_\Gamma  (I) \le \Gamma $ of elements $\gamma  \in \Gamma $ such that $I = I^\gamma  $
     is said to be the normalizer of the ideal $I$ in the group $\Gamma $.   The subgroup
     $S_\Gamma  (I) \le \Gamma $ of elements $\gamma  \in \Gamma $ such that
     $I \cap kB = I^\gamma   \cap kB$ for some finitely generated dense subgroup $B$ of $A$ is said
     to be the standardiser of the ideal $I$ in the group $\Gamma $. Evidently, $N_\Gamma  (I) \le S_\Gamma  (I)$.
     Theorem 4.4 states that if the field $k$ has characteristic zero and $P$ is a prime faithful ideal of
     the group algebra $kA$ such that $S_\Gamma  (P) = \Gamma $ then $P$ is controlled by $\Delta _\Gamma  (A)$.
     In the case, where the group $A$ is finitely generated and $\left| {\Gamma :N_\Gamma  (P)} \right| < \infty $,
     such a result was proved by Roseblade in [7, theorem D]. As lemma 4.2 shows, it is sufficient to consider
     the situation where the field $k$ is finitely generated, the group $A$ is minimax and the group $\Gamma $
     is cyclic. Lemma 4.3 shows that it is sufficient to prove that the ideal $P$ is controlled by a finitely
     generated $\Gamma $-invariant subgroup of $A$. So, if the ideal $P$ is $\Gamma $-invariant, the result
     would follow immediately from theorem 3.1 or from [8, theorem 1.1], because the controller of a
     $\Gamma $-invariant ideal is a $\Gamma $-invariant subgroup. However, the condition
     $S_\Gamma  (P) = \Gamma $ is much more general than $N_\Gamma  (P) = \Gamma $ and
     does not mean that the ideal $P$ is $\Gamma $-invariant. This circumstance creates
     the main difficulties which are conquered in the proof of theorem 4.4.
\vspace{1 pc}\par
     {\bf Lemma 4.1.} {\it Let $k$ be a field, let $A$ be a torsion-free abelian group of finite rank
     acted by a group $\Gamma$ and let $P$ be a faithful prime ideal of the group algebra $kA$.
     Then:\\
     (i) $r(c(P)) = r(c(kB \cap P))$ for any dense subgroup $B$ of $A$;\\
     (ii) if  $S_\Gamma  (P) = \Gamma $ then $is_A c(P)$ is a $\Gamma $-invariant subgroup of  $A$.}\par
     {\bf Proof.} (i) Since $c(kB \cap P) \le c(P) \cap B$, there is no harm in assuming that $is_A c(P) = A$.
     Suppose that $r(c(P)) > r(c(kB \cap P))$ then there is an isolated subgroup $D$ of $A$ such that
     $c(P \cap kB) \le D \cap B$ and $r(A/D) = 1$. Let $C$ be a maximal isolated subgroup of $D$ such
     that $P \cap kC = 0$.  Since $c(P)$ is a dense subgroup of $A$, we can conclude that $C$ is a
     maximal subgroup of $A$ such that  $kC \cap P = 0$. Then the transcendence degree over $k$ of
     the field of fraction $k_1 $ of the domain $kA/P$ is $r(C)$. On the other hand, as
     $k(C \cap B) \cap P = 0$ and $c(P \cap kB) \le D \cap B$, we see that the transcendence
     degree of the field of fractions $k_2 $ of the domain $kB/(kB \cap P)$ is at least $r(C) + 1$.
     But it is impossible because $A/B$ is a torsion group and hence $k_1 $ is an algebraic
     extension of $k_2 $.\par
     (ii) By the definition of $S_\Gamma  (P)$, for any element $\gamma  \in \Gamma $ there is
     a finitely generated dense subgroup $B$ of $A$ such that $kB \cap P = kB \cap P^\gamma  $.
     Therefore $c(kB \cap P) \le c(P) \cap c(P)^\gamma  $, and by (i),
     $r(c(P)) = r(c(kB \cap P)) = r(c(P)^\gamma  )$. It easily implies that $c(P) \cap c(P)^\gamma
     $ is a dense subgroup in $c(P)$ and $c(P)^\gamma  $. Therefore,   $is_A c(P) = is_A c(P)^\gamma  $.
\vspace{1 pc}\par
     {\bf Lemma 4.2.}{\it Suppose that there exist a field $k$, a
     torsion-free abelian group $A$ of finite rank acted by a group $\Gamma $
     and a faithful prime ideal $P$   of $kA$ which is not controlled by
     $\Delta _\Gamma  (A)$ and such  that $S_\Gamma  (P) = \Gamma $. Then
     there exist a finitely generated subfield $k_1  \le k$, a minimax subgroup
     $A_1  \le A$ acted by a cyclic subgroup $\left\langle \gamma  \right\rangle  \le \Gamma $
      and a faithful prime ideal $P_1  = P \cap k_1 A_1 $   of  $k_1 A_1 $ which is not
      controlled by  $\Delta _{\left\langle \gamma  \right\rangle } (A_1 )$ and such
      that $S_{\left\langle \gamma  \right\rangle } (P) = \left\langle \gamma  \right\rangle
      $.}\par
     {\bf Proof.} At first, we show that there is no harm in assuming that the field $
     K$ is finitely generated. Suppose that $P$ is not controlled by $\Delta _\Gamma  (A)$
     then there is an element $\alpha  \in P\backslash ((P \cap k\Delta _\Gamma  (A))kA){\rm{ }}$.
     Let $k_1 $ be a subfield of  $k$ generated by coefficients of $\alpha $ and let $P_1  = P \cap k_1 A$
     then $\alpha  \in P_1 \backslash ((P_1  \cap k_1 \Delta _\Gamma  (A))k_1 A){\rm{ }}$ and hence $P_1 $
     is not controlled by $\Delta _\Gamma  (A)$. It is easy to note that $P_1 $ is a prime faithful ideal
     of $k_1 A$ and $S_\Gamma  (P_1 ) = \Gamma $. Thus, it is sufficient to consider the case where the
     field  $k$ is finitely generated.\par
     Let $c(P)$ be the controller of $P$, by lemma 4.1(ii), $is_A c(P)$ is a $\Gamma $-invariant subgroup
     of $A$. Then there is no harm in assuming that $is_A c(P) = A$. We also assume that $\Gamma $ acts
     on the group $A$ faithfully, that is $C_\Gamma  (A) = 1$.\par
     Suppose that $c(P)$ is not contained in $\Delta _\Gamma  (A)$. Let $F$ be a free dense subgroup of
     $c(P)$ with free generators $\left\{ {a_i |i = 1,...,n}
     \right\}$.\par
     Suppose that for any element $\gamma  \in \Gamma $ we have
     $\left| {\left\langle \gamma  \right\rangle /C_{\left\langle \gamma  \right\rangle } (a_i )} \right| < \infty $.
      Let $C_\gamma   =  \cap _{i = 1}^n C_{\left\langle \gamma  \right\rangle } (a_i )$ then
      $\left| {\left\langle \gamma  \right\rangle /C_\gamma  } \right| < \infty $ and it easy
      to see that $C_\gamma  $ centralizes $A$. Since $\Gamma $ acts on the group $A$ faithfully,
      it implies that the group $\Gamma $  is torsion. But $\Gamma$ is a linear group over the
      field of rational numbers and hence the group $\Gamma $ is finite. Then, evidently,
      $\Delta _\Gamma  (A) = A$ and a contradiction is
      obtained.\par
      Thus, if $c(P)$ is not
      contained in $\Delta _\Gamma  (A)$, then there are an element $a \in c(P)$ and an element
      $\gamma  \in \Gamma $ such that the group $\left\langle \gamma  \right\rangle $ is infinite
      cyclic and $C_{\left\langle \gamma  \right\rangle } (a) = 1$. Then we can replace $\Gamma $
      by $\left\langle \gamma  \right\rangle $ because $c(P)$ is not contained in
      $\Delta _{\left\langle \gamma  \right\rangle } (A)$. Let $A_1 $ be a dense subgroup of $A$
      which is finitely generated as a $\left\langle \gamma  \right\rangle $-module and which
      contains the element $a$. Then, by [3, lemma 5.2], the subgroup $A_1 $ is
      minimax.\par
      Put $P_1  = P \cap kA_1 $,  by lemma 4.1(i), $c(P_1 )$ is a dense subgroup of$A_1 $. If $c(P_1 )$
      is contained in $\Delta _\Gamma  (A_1 )$ then $\Delta _\Gamma  (A_1 )$ is a dense subgroup of
      $A_1 $ and  hence $\Delta _\Gamma  (A_1 ) = A_1 $ . But it is impossible because
      $a \notin \Delta _\Gamma  (A_1 )$. Thus, $c(P_1 )$ is not contained in $\Delta _\Gamma  (A_1 )$.
\vspace{1 pc}\par

     {\bf Lemma 4.3.}{\it Let $k$ be a field of characteristic zero, let $A$ be a torsion-free finitely generated
     abelian group acted by a cyclic group $\Gamma  = \left\langle g \right\rangle $ and let $P$ be a
     faithful prime ideal of the group algebra $kA$. Suppose that $S_\Gamma  (P) = \Gamma $ then $P$
     is controlled by $\Delta _\Gamma  (A)$.}\par
     {\bf Proof.}  Since $S_\Gamma  (P) = \Gamma $, there is a subgroup $B$ of finite index in $A$ such that
     $P \cap kB = P^\gamma   \cap kB$. As $A^n  \le B$ fore some positive integer $n$, we can assume
     that the subgroup $B$ is $\Gamma $-invariant and it implies that the ideal $P_1  = P \cap kB$ is
     $\left\langle \gamma  \right\rangle $-invariant. Since $\left| {A:B} \right| < \infty $, there is
     only finite set of ideals $X$ of $kA$ such that $P_1  = X \cap kB$. Then, as the ideal $P_1 $ is
     $\left\langle \gamma  \right\rangle $-invariant and $\Gamma  = \left\langle g \right\rangle $,
     it easily implies that $\left| {\Gamma :N_\Gamma  (P)} \right| < \infty $   and the assertion
     follows from [7, theorem D].
\vspace{1 pc}\par

     {\bf Theorem 4.4.}{\it Let $k$ be a field of characteristic zero, let $A$ be a torsion-free abelian group of
     finite rank acted by a group $\Gamma $ and let $P$ be a faithful prime ideal of the group algebra
     $kA$. Suppose that $S_\Gamma  (P) = \Gamma $ then  $P$ is controlled by $\Delta _\Gamma
     (A)$.}\par
     {\bf Proof.} By lemma 4.2, it is sufficient to consider the case where the field $k$ is finitely
     generated, the group $A$ is minimax and the group $\Gamma $ is
     cyclic.\par
     By lemma 4.1(i), $is_A c(P)$ is a $\Gamma $-invariant subgroup of  $A$ and changing $A$
     by $is_A c(P)$, we can assume that $c(P)$ is a dense subgroup of $A$. Then it follows
     from theorem 3.1 that $c(P)$ is a finitely generated dense subgroup of $A$.\par
     Since the group $A$ is minimax, the set $\pi (A/c(P))$ is finite. Let $\tau $ be the set of
     all roots from 1 of degree $p^n $ for all $p \in \pi (A/c(P))$ and all positive integer $n$
     in addition with a primitive root from 1 of degree 4.\par
     Put $C = c(P)$. Since $S_\Gamma  (P) = \Gamma $,
     for any $\gamma  \in \Gamma $ there is a finitely generated dense subgroup $D_\gamma  $ of  $A$
     such that $kD_\gamma   \cap P = kD_\gamma   \cap P^\gamma  $. Evidently, we can assume that
     $D_\gamma   \le C \cap C^\gamma  $.\par
     Let $R = kA/P$ and let $h$ be the field of fractions of $R$, then $A \le h^ *  $. Let $f$ be
     a field obtained by addition to the field $h$ of all roots from the set $\tau $ and let $s$
     be a subfield of $f$ generated by all roots from the set $\tau $ and the subgroup $C$. Let
     $B = A \cap s^ *  $ then, by theorem 3.3, $B$ is a finitely generated dense subgroup of $A$
     which controls $P$. Let $d$ be a subfield of $f$ generated by $D_\gamma  $ and $t(s^*)$ and
     let $H = A \cap d^ *  $ then, by theorem 3.3, $c(P) \le H \le  B$.\par
     Let $R_\gamma   = kA/P^\gamma  $ and let $h_\gamma  $ be the field of fractions of $R_\gamma  $.
     Let $f_\gamma  $ be a field obtained by addition to the field $h_\gamma  $ all roots from the
     set $\tau $ and let $s_\gamma  $ be a subfield of $f_\gamma  $ generated by all roots from the
     set $\tau $ and the subgroup $C^\gamma  $.  Let $d_\gamma  $ be a subfield of $f_\gamma  $
     generated by $D_\gamma  $ and $t(s_\gamma  ^*)$ and let $H_\gamma   = A \cap (d_\gamma^ *)$ then,
     by theorem 3.3, $C^\gamma   = c(P^\gamma  ) \le H_\gamma $.\par
     Since $(kA/P)^\gamma   = kA/P^\gamma  $, we can conclude that the fields $f$ and $f_\gamma  $ are
     isomorphic under the isomorphism $\varphi $ induced by action of $\gamma $. Since
     $\varphi (C) = C^\gamma  $, we see that the fields $s$ and $s_\gamma  $ are isomorphic
     and hence  $t(s^*) \cong t(s_\gamma  ^*)$.\par
     Let $K_1 $ be a subring of $f$ generated by $k$ and $D_\gamma  $,  and let $K_2 $ be a
     subring of $f_\gamma  $ generated by $k$ and $D_\gamma  $. Since
     $kD_\gamma   \cap P = kD_\gamma   \cap P^\gamma  $, we can conclude that $K_1  = K_2  = K$.
     Then, as $t(s^*) \cong t(s_\gamma  ^*)$,  the fields $d$ and $d_\gamma  $ are the fields of
     decomposition of the same set $\Omega $ of polynomials over the field of fractions of the
     domain $K$.  More precisely, $\Omega $ is the set of  cyclotomic polynomials for roots from 1
     which belong to the set  $t(s^*) \cong t(s_\gamma  ^*)$. Therefore, by [6, Chap. VII, theorem 3],
     there is an isomorphism $\psi $ between the fields $d$ and $d_\gamma  $ which centralizes the
     elements of $K$ and we can assume that $K \le d \cap d_\gamma $.\par
     Put $\pi  = \pi (A/C)$ then $\pi  = \pi (A/C^\gamma  )$. As $C \le H$ and $C^\gamma   \le H_\gamma  $,
     we can conclude that $A/H$ and $A/H_\gamma  $ are $\pi     $-groups.\par
     Suppose that $H \ne H_\gamma  $
     then either $H\backslash H_\gamma   \ne \emptyset $ or $H_\gamma  \backslash H \ne \emptyset $.
     Suppose that $H\backslash H_\gamma   \ne \emptyset $ then $H/(H \cap H_\gamma  )$ is a nontrivial
     $\pi $-group and hence the Hall $\pi $-subgroup of the quotient group $H/D_\gamma  $ is not contained in
     $(H \cap H_\gamma  )/D_\gamma  $. Therefore, there exists an element $a \in H\backslash H_\gamma  $
     such that $\left| {\bar a} \right| = m$ is a $\pi $-number, that is all prime divisors of $m$ belong
     to $\pi $, where $\bar a$ is the image of $a$ in the quotient group $H/D_\gamma  $. Thus, $m$ is the
     smallest integer such that $a^m  = b \in D_\gamma  $. Since $a \in H = A \cap d$ and $d$ contains all
     roots from the set $\tau $, we can conclude that all roots of the polynomial $X^m  - b$ are contained
     in $d$. On the other hand, since $a \notin H_\gamma   = A \cap d_\gamma  $ and $a$ is a root of
     polynomial $X^m  - b$ in $f_\gamma  $, we can conclude that $d_\gamma  $ does not contain all
     roots of the polynomial $X^m  - b$. Since the fields $d$ and $d_\gamma  $ are isomorphic and
     $b \in K \le d \cap d_\gamma  $, it leads to a contradiction. If $H_\gamma  \backslash H \ne \emptyset $
     then the same arguments lead to the same contradiction. Thus, $H = H_\gamma
     $.\par
     Evidently, $c(P^\gamma  ) = (c(P))^\gamma  $. It follows from theorem 3.3 that
     $(c(P))^\gamma   \le H_\gamma  $ and, as $H_\gamma   = H \le B$, we can conclude that
     $(c(P))^\gamma   \le B$ for each $\gamma  \in \Gamma $. Then
     $S = \left\langle {{(c(P))^\gamma  }} \mathrel{\left | {\vphantom {{(c(P))^\gamma  }
     {\gamma  \in \Gamma }}} \right. \kern-\nulldelimiterspace} {{\gamma  \in \Gamma }} \right\rangle  \le B$
     and hence $S$ is a finitely generated $\Gamma $-invariant subgroup of $A$ which controls $P$. Then
     the assertion follows from lemma 4.3.
\vspace{1 pc}\par
\begin{center}
     REFERENCES
\end{center}
\vspace{1 pc}
 1. Brookes, Ch.J.B. Ideals in group rings of
soluble groups of finite rank. Math. Proc. Camb. Phil. Soc. 1985,
V. 97,  P. 27-49.\\
2. Brown, K.S. Cohomology of groups. New York-Berlin:
Springer-Verlag; 1982.\\
 3. Hall, P. On the finiteness of certain
soluble groups. Proc. London Math. Soc. 1959, V.9, No 36, P. 595 -
622.\\
 4. Janusz, G.J. Algebraic number fields. Providance: AMS;
1996. \\
5. Lambek, J. Rings and Modules. Waltham – Massachusetts – Toronto
– London: Blaisdell publishing company; 1966.\\
 6. Lang, S.
Algebra. Mass.: Addison-Wesley Publ. Comp. Reading; 1965.\\
 7.Roseblade, J.E. Prime ideals in group rings of polycyclic groups.
Proc. London Math. Soc. 1976, V. 36, No. 3, P.385-447. \\
8. Segal,D.  On the group rings of abelian minimax groups. J.
Algebra, 2001,- V.237,  P.64-94. \\
9. Tushev, A.V. Noetherian modules over abelian groups of finite
torsion-free rank. Ukrainian Math. J. 1991, V.43, No 7,8 , P.
975-981.\\
 10. Tushev, A.V. Noetherian modules over minimax
abelian groups. Ukrainian Math. Zh. 2002, V. 54, No 7, P. 974-985.
(in Russian) \\
11. Tushev, A.V. Induced representations of abelian groups of
finite rank. Ukrainian Math. Zh. 2003, V. 55, No 9, P. 974-985.
(in Russian) \\
12. Tushev, A.V. On deviation in groups. Illinois
Journal of Mathematics, 2003, V. 47, No 1/2, P. 539-550.

\end{document}